\pdfoutput=1
\RequirePackage{ifpdf}
\ifpdf 
\documentclass[pdftex]{sigma}
\else
\documentclass{sigma}
\fi

\newcommand{\eqdef}{\stackrel{\text{def}}{=}}
\newcommand{\Romannumeral}[1]{\uppercase\expandafter{\romannumeral#1}}
\newcommand{\I}{\text{\Romannumeral{1}}}
\newcommand{\II}{\text{\Romannumeral{2}}}

\newcommand{\cF}{c_{\text{\tiny$\mathcal{F}$}}}

\begin{document}

\newcommand{\arXivNumber}{1612.00927}

\renewcommand{\PaperNumber}{020}

\FirstPageHeading

\ShortArticleName{Simplif\/ied Expressions of the Multi-Indexed Laguerre and Jacobi Polynomials}

\ArticleName{Simplif\/ied Expressions of the Multi-Indexed Laguerre and Jacobi Polynomials}

\Author{Satoru ODAKE and Ryu SASAKI}
\AuthorNameForHeading{S.~Odake and R.~Sasaki}
\Address{Faculty of Science, Shinshu University, Matsumoto 390-8621, Japan}
\Email{\href{mailto:odake@azusa.shinshu-u.ac.jp}{odake@azusa.shinshu-u.ac.jp}, \href{mailto:ryu@yukawa.kyoto-u.ac.jp}{ryu@yukawa.kyoto-u.ac.jp}}

\ArticleDates{Received December 30, 2016, in f\/inal form March 23, 2017; Published online March 29, 2017}

\Abstract{The multi-indexed Laguerre and Jacobi polynomials form a complete set of orthogonal polynomials. They satisfy second-order dif\/ferential equations but not three term recurrence relations, because of the `holes' in their degrees. The multi-indexed Laguerre and Jacobi polynomials have Wronskian expressions originating from multiple Darboux transformations. For the ease of applications, two dif\/ferent forms of simplif\/ied expressions of the multi-indexed Laguerre and Jacobi polynomials are derived based on various identities. The parity transformation property of the multi-indexed Jacobi polynomials is derived based on that of the Jacobi polynomial.}


\Classification{42C05; 33C45; 34A05}

\section{Introduction}\label{sec:intro}

The exceptional and multi-indexed orthogonal polynomials \cite{quesne+,gomez+,gomez,gomez2,os16,os25,os26,os27,os35,quesne} seem to be a focal point of recent research on exactly solvable quantum mechanics. They belong to a new type of orthogonal polynomials which satisfy second-order dif\/ferential (dif\/ference) equations and form complete orthogonal basis in an appropriate Hilbert space. One of their characteristic features is that they do not satisfy three term recurrence relations because of the `holes' in the degrees. This is how they avoid the constraints due to Bochner \cite{bochner, routh}. They are constructed from the original quantum mechanical systems, the radial oscillator potential and the P\"oschl--Teller potential, by multiple application of Darboux transformations \cite{crum, darboux} in terms of seed solutions called virtual state wavefunctions which are generated by two types of discrete symmetry transformations \cite{os25,os26,os27,os35}.

The multi-indexed Laguerre and Jacobi polynomials have Wronskian expressions \cite{os25} originating from multiple Darboux transformations~\cite{crum}. In this note we present two dif\/ferent forms of equivalent determinant expressions without higher derivatives of the Wronskians by using various identities of the Laguerre and Jacobi polynomials \cite{szego}. These simplif\/ied expressions show explicitly the constituents of the multi-indexed orthogonal polynomials and they are helpful for deeper understanding. In \cite{d14} and \cite{d17} Dur\'an employed similar simplif\/ied expressions as the starting point of his exposition of the exceptional Laguerre and Jacobi polynomials. See also \cite{dp15} by Dur\'an and P\'erez. In their expressions, the matrix elements of the determinants are polynomials. In the original expressions in \cite{os25}, the matrix elements of the determinants contain the non-polynomial factors, see \eqref{XiD}--\eqref{Pn}. In the simplif\/ied expressions presented in this paper in Sections~\ref{sec:simpA} and~\ref{sec:simpB}, they are also all polynomials.

This short note is organised as follows. In Section~\ref{sec:ori} the quantum mechanical settings for the original Laguerre and Jacobi polynomials are recapitulated. That is, the Hamiltonians with the radial oscillator potential and the P\"oschl--Teller potential are introduced and their eigenvalues and eigenfunctions are presented. Type~$\I$ and~$\II$ discrete symmetry transformations for the Hamiltonians of the Laguerre and Jacobi polynomials are explained in Section~\ref{sec:sym}. The seed solutions for Darboux transformations, to be called the virtual state wavefunctions of Type~$\I$ and~$\II$, for the Laguerre and Jacobi, are listed explicitly. In Section~\ref{sec:wron} the Wronskian forms of the multi-indexed Laguerre and Jacobi polynomials derived in \cite{os25} and \cite{rrmiop3} are recapitulated as the starting point. Sections~\ref{sec:simpA} and~\ref{sec:simpB} are the main content of this note. Those who are familiar with the multi-indexed Laguerre and Jacobi polynomials can directly go to this part. The f\/irst simplif\/ied expressions of the multi-indexed Laguerre and Jacobi polynomials are derived in Section~\ref{sec:simpA} by using various identities of the original Laguerre and Jacobi polynomials. The second simplif\/ied expressions, to be derived in Section~\ref{sec:simpB}, are the consequences of the multi-linearity of determinants and the form of the Schr\"odinger equation \smash{$\psi''(x)=(U(x)-\mathcal{E})\psi(x)$}. Every even order derivative $\psi^{(2m)}(x)$ in the Wronskian can be replaced by $(-\mathcal{E})^m\psi(x)$. In Section~\ref{sec:pari} the parity transformation formula of the multi-indexed Jacobi polynomials is presented. Section~\ref{sec:summary} is for a summary and comments.

\section{Multi-indexed Laguerre and Jacobi polynomials}\label{sec:MILJ}

The foundation of the theory of multi-indexed orthogonal polynomials is the exactly solvable one-dimensional quantum mechanical system $\mathcal{H}$:
\begin{gather}
 \mathcal{H}\phi_n(x)=\mathcal{E}_n\phi_n(x), \quad n=0,1,\ldots,\qquad \mathcal{H}=-\frac{d^2}{dx^2}+U(x), \label{scheq}
\end{gather}
and its iso-spectrally deformed system
\begin{gather}
 \mathcal{H}_{\mathcal{D}}\phi_{\mathcal{D}\, n}(x) =\mathcal{E}_n\phi_{\mathcal{D}\,n}(x), \quad n=0,1,\ldots,\qquad
 \mathcal{H}_{\mathcal{D}}=-\frac{d^2}{dx^2}+U_{\mathcal{D}}(x), \label{Dscheq}
\end{gather}
in terms of multiple application of Darboux transformations~\cite{crum, os25}. In this note we discuss two systems which have the Laguerre and Jacobi polynomials as the main parts of the eigenfunctions. We follow the notation of~\cite{os25} with slight modif\/ications for simplif\/ication sake.

\subsection{Original Laguerre and Jacobi polynomials}\label{sec:ori}

\subsubsection{Radial oscillator potential}

The Hamiltonian with the radial oscillator potential
\begin{gather}
 U(x)\eqdef x^2+\frac{g(g-1)}{x^2}-2g-1,\qquad 0<x<\infty,\qquad g>\frac12, \label{rad}
\end{gather}
has the Laguerre polynomials $L^{(\alpha)}_n(\eta)$ as the main part of the eigenfunctions:
\begin{gather*}
 \phi_n(x;g)\eqdef\phi_0(x;g)L^{(g-\frac12)}_n\bigl(\eta(x)\bigr),\qquad \mathcal{E}_n\eqdef 4n,\qquad\eta(x)\eqdef x^2,\\
 \phi_0(x;g)\eqdef e^{-\frac12x^2}x^g,\qquad L_n^{(\alpha)}(\eta)=\frac{1}{n!} \sum_{k=0}^n\frac{(-n)_k}{k!}(\alpha+k+1)_{n-k}\eta^k.
\end{gather*}
Here are some identities of the Laguerre polynomials \cite{szego} to be used in Sections~\ref{sec:simpA} and \ref{sec:simpB},
\begin{gather}
\partial_{\eta}L_n^{(\alpha)}(\eta)=-L_{n-1}^{(\alpha+1)}(\eta), \label{Lforward}\\
L_n^{(\alpha)}(\eta)-L_n^{(\alpha-1)}(\eta)=L_{n-1}^{(\alpha)}(\eta), \label{Lid1}\\
\eta L_{n-1}^{(\alpha+1)}(\eta)-\alpha L_{n-1}^{(\alpha)}(\eta) =-nL_n^{(\alpha-1)}(\eta). \label{Lid2}
\end{gather}
As is clear from \eqref{rad}, the lower boundary, $x=0$, is the regular singular point with the characteristic exponents $g$ and $1-g$. The upper boundary point, $x=\infty$, is an irregular singular point. Here $(a)_n\eqdef\prod\limits_{j=1}^n(a+j-1)$ is the shifted factorial or the so-called Pochhammer's symbol.

\subsubsection{P\"oschl--Teller potential}

The Hamiltonian with the P\"oschl--Teller potential
\begin{gather*}
 U(x)\eqdef\frac{g(g-1)}{\sin^2x}+\frac{h(h-1)}{\cos^2x}-(g+h)^2,\qquad 0<x<\frac{\pi}{2},\qquad g>\frac12,\qquad h>\frac12,
\end{gather*}
has the Jacobi polynomials $P_n^{(\alpha,\beta)}(\eta)$ as the main part of the eigenfunctions:
\begin{gather*}
\phi_n\bigl(x;(g,h)\bigr)\eqdef\phi_0\bigl(x;(g,h)\bigr) P_n^{(g-\frac12,h-\frac12)}\bigl(\eta(x)\bigr),\qquad \eta(x)\eqdef\cos2x,\\
\phi_0\bigl(x;(g,h)\bigr)\eqdef(\sin x)^g(\cos x)^h,\qquad \mathcal{E}_n(g,h)\eqdef 4n(n+g+h),\\
P_n^{(\alpha,\beta)}(\eta)=\frac{(\alpha+1)_n}{n!} \sum_{k=0}^n\frac{1}{k!}\frac{(-n)_k(n+\alpha+\beta+1)_k}{(\alpha+1)_k} \left(\frac{1-\eta}{2}\right)^k.
\end{gather*}
Here are some identities of the Jacobi polynomials \cite{szego} to be used in in Sections~\ref{sec:simpA} and~\ref{sec:simpB},
\begin{gather}
\partial_{\eta}P_n^{(\alpha,\beta)}(\eta) =\tfrac12(n+\alpha+\beta+1)P_{n-1}^{(\alpha+1,\beta+1)}(\eta), \label{Jforward}\\
(n+\beta)P_n^{(\alpha,\beta-1)}(\eta) -\beta P_n^{(\alpha-1,\beta)}(\eta)
 =(n+\alpha+\beta)\tfrac{1+\eta}{2}P_{n-1}^{(\alpha,\beta+1)}(\eta), \label{Jid1}\\
(n+\alpha)P_n^{(\alpha-1,\beta)}(\eta) -\alpha P_n^{(\alpha,\beta-1)}(\eta)
 =-(n+\alpha+\beta)\tfrac{1-\eta}{2}P_{n-1}^{(\alpha+1,\beta)}(\eta). \label{Jid1m}
\end{gather}
The two boundary points, $x=0,\frac{\pi}2$ are the regular singular points with the characteristic exponents $g$, $1-g$ and $h$, $1-h$, respectively.

\subsection{Discrete symmetry transformations and virtual state wavefunctions}\label{sec:sym}

The seed solutions (virtual state wavefunctions) for Darboux transformations can be constructed by applying discrete symmetry transformations of the
Hamiltonian to the eigenfunctions. They have negative energies and have no node and they are not square integrable, see~\cite{os25} for more detail.

\subsubsection{Laguerre (L) system}

{\bf Type $\I$ transformation.} It is obvious that the transformation $x\to ix$ is the symmetry of the radial oscillator system. The seed solutions are
\begin{alignat*}{3}
 & \text{L1}\colon \quad && \mathcal{H}\tilde{\phi}^{\I}_{\text{v}}(x) =\tilde{\mathcal{E}}^{\I}_{\text{v}}\tilde{\phi}^{\I}_{\text{v}}(x),\qquad
 \tilde{\phi}^{\I}_{\text{v}}(x)\eqdef \tilde{\phi}^{\I}_0(x)\xi^{\I}_{\text{v}}\bigl(\eta(x)\bigr),\qquad
 \tilde{\phi}^{\I}_0(x)\eqdef e^{\frac12x^2}x^g,&\nonumber\\
&&& \xi^{\I}_{\text{v}}(\eta)\eqdef L^{(g-\frac12)}_{\text{v}}(-\eta), \qquad\tilde{\mathcal{E}}^{\I}_{\text{v}}\eqdef-4(g+\text{v}+\tfrac12),\qquad \text{v}\in\mathbb{Z}_{\ge0}.&
\end{alignat*}
Hereafter we use $\text{v}$ for the degree of the seed polynomial $\xi_{\text{v}}$ in order to distinguish it with the degree~$n$ of eigenpolynomial~$P_n$.

\noindent
{\bf Type $\II$ transformation.} The exchange of the characteristic exponents of the regular singular point $x=0$, $g\to1-g$, is another discrete symmetry transformation. It generates the seed solutions,
\begin{alignat*}{3}
& \text{L2}\colon \quad && \mathcal{H}\tilde{\phi}^{\II}_{\text{v}}(x)
 =\tilde{\mathcal{E}}^{\II}_{\text{v}}\tilde{\phi}^{\II}_{\text{v}}(x),\qquad
 \tilde{\phi}^{\II}_{\text{v}}(x)\eqdef \tilde{\phi}^{\II}_0(x)\xi^{\II}_{\text{v}}\bigl(\eta(x)\bigr),\qquad
 \tilde{\phi}^{\II}_0(x)\eqdef e^{-\frac12x^2}x^{1-g},&\nonumber\\
&&& \xi^{\II}_{\text{v}}(\eta)\eqdef L^{(\frac12-g)}_{\text{v}}(\eta),
 \qquad\tilde{\mathcal{E}}^{\II}_{\text{v}}\eqdef-4\big(g-\text{v}-\tfrac12\big),\qquad
 \text{v}=0,1,\ldots,\big[g-\tfrac12\big]',&
\end{alignat*}
in which $[a]'$ denotes the greatest integer less than $a$.

\subsubsection{Jacobi (J) system}

{\bf Type $\I$ transformation.} The exchange of the characteristic exponents of the regular singular point $x=\frac{\pi}{2}$, $h\to1-h$, is a~discrete symmetry transformation:
\begin{alignat*}{3}
& \text{J1} \colon \quad && \mathcal{H}\tilde{\phi}^{\I}_{\text{v}}(x)
 =\tilde{\mathcal{E}}^{\I}_{\text{v}}\tilde{\phi}^{\I}_{\text{v}}(x),\qquad
 \tilde{\phi}^{\I}_{\text{v}}(x)\eqdef \tilde{\phi}^{\I}_0(x)\xi^{\I}_{\text{v}}\bigl(\eta(x)\bigr),\qquad
 \tilde{\phi}^{\I}_0(x)\eqdef(\sin x)^g(\cos x)^{1-h},& \nonumber\\
&&& \xi^{\I}_{\text{v}}(\eta)\eqdef P^{(g-\frac12,\frac12-h)}_{\text{v}}(\eta), \qquad \tilde{\mathcal{E}}^{\I}_{\text{v}}\eqdef
 -4\big(g+\text{v}+\tfrac12\big)\big(h-\text{v}-\tfrac12\big),& \\
 &&& \text{v}=0,1,\ldots,\big[h-\tfrac12\big]'.&
\end{alignat*}

\noindent
{\bf Type $\II$ transformation.} Likewise the exchange of the characteristic exponents of the regular singular point $x=0$, $g\to1-g$, is a discrete symmetry transformation:
\begin{alignat*}{3}
 & \text{J2}\colon \quad && \mathcal{H}\tilde{\phi}^{\II}_{\text{v}}(x) =\tilde{\mathcal{E}}^{\II}_{\text{v}}\tilde{\phi}^{\II}_{\text{v}}(x),\qquad
 \tilde{\phi}^{\II}_{\text{v}}(x)\eqdef \tilde{\phi}^{\II}_0(x)\xi^{\II}_{\text{v}}\bigl(\eta(x)\bigr),\qquad
 \tilde{\phi}^{\II}_0(x)\eqdef(\sin x)^{1-g}(\cos x)^h,&\nonumber\\
&&& \xi^{\II}_{\text{v}}(\eta)\eqdef P^{(\frac12-g,h-\frac12)}_{\text{v}}(\eta), \qquad \tilde{\mathcal{E}}^{\II}_{\text{v}}\eqdef
 -4\big(g-\text{v}-\tfrac12\big)\big(h+\text{v}+\tfrac12\big),&\\
 &&& \text{v}=0,1,\ldots,\big[g-\tfrac12\big]'.&
\end{alignat*}

\subsection{Wronskian forms of the multi-indexed Laguerre and Jacobi polynomials}\label{sec:wron}

We deform the original Hamiltonian system \eqref{scheq} by applying multiple Darboux transforma\-tions~\cite{crum, os25} in terms of $M_{\I}$ Type~$\I$ seed solutions and $M_{\II}$ Type~$\II$ seed solutions specif\/ied by the degree set $\mathcal{D}\eqdef\{d_1,\ldots,d_M\}$ (ordered set), $M\eqdef M_{\I}+M_{\II}$, $M_{\I}\eqdef\#\{d_j\,|\,d_j\in\mathcal{D},\,\text{Type $\I$}\}$, $M_{\II}\eqdef\#\{d_j\,|\,d_j\in\mathcal{D},\,\text{Type $\II$}\}$. The multi-indexed orthogonal polynomials $\{P_{\mathcal{D},n}(\eta)\}$ are the main parts of the eigenfunctions $\{\phi_{\mathcal{D}\,n}(x)\}$ of the deformed Hamiltonian system $\mathcal{H}_{\mathcal{D}}$~\eqref{Dscheq}:
\begin{gather*}
 U_{\mathcal{D}}(x)\eqdef U(x)-2\partial_x^2 \log\bigl|\text{W}\big[\tilde{\phi}_{d_1},\ldots,\tilde{\phi}_{d_M}\big](x)\bigr|,\\
 \phi_{\mathcal{D}\,n}(x)\eqdef \frac{\text{W}\big[\tilde{\phi}_{d_1},\ldots,\tilde{\phi}_{d_M},\phi_n\big](x)}
 {\text{W}\big[\tilde{\phi}_{d_1},\ldots,\tilde{\phi}_{d_M}\big](x)}\eqdef\cF^{M}\psi_{\mathcal{D}}(x)
 P_{\mathcal{D},n}\bigl(\eta(x)\bigr),\qquad n=0,1,\ldots,\\
 \psi_{\mathcal{D}}(x)\eqdef\frac{\hat{\phi}_0(x)} {\Xi_{\mathcal{D}}\bigl(\eta(x)\bigr)},
\end{gather*}
in which $\hat{\phi}_0(x)$ and $\cF$ are def\/ined by
\begin{gather*}
 \hat{\phi}_0(x)\eqdef
 \begin{cases}
 \phi_0(x;g+M_{\I}-M_{\II}),&\text{L},\\
 \phi_0\bigl(x;(g+M_{\I}-M_{\II},h-M_{\I}+M_{\II})\bigr),& \text{J},
 \end{cases}\qquad
 \cF\eqdef
 \begin{cases}
 \hphantom{-}2,& \text{L},\\
 -4, & \text{J}.
 \end{cases}
\end{gather*}
The superscript $\I$ and $\II$ of the seed solutions are suppressed for simplicity of notation. The Wronskian of $n$-functions $\{f_1,\ldots,f_n\}$ is def\/ined by formula
\begin{gather*}
 \text{W}[f_1,\ldots,f_n](x) \eqdef\det\left(\frac{d^{j-1}f_k(x)}{dx^{j-1}}\right)_{1\leq j,k\leq n}.
\end{gather*}

{\allowdisplaybreaks There are two dif\/ferent but equivalent Wronskian def\/initions of the denominator polynomial $\Xi_{\mathcal{D}}(\eta)$ and the multi-indexed orthogonal polynomial $P_{\mathcal{D},n}(\eta)$. The f\/irst is based on the Wronskians of the `polynomials' \cite{os25}:
 \begin{gather}
 \Xi_{\mathcal{D}}(\eta) \eqdef\text{W}[\mu_{d_1},\ldots,\mu_{d_M}](\eta)\times
 \begin{cases}
 \eta^{(M_{\I}+g-\frac12)M_{\II}}e^{-M_{\I}\eta},&\text{L},\\
 \bigl(\frac{1-\eta}{2}\bigr)^{(M_{\I}+g-\frac12)M_{\II}}
 \bigl(\frac{1+\eta}{2}\bigr)^{(M_{\II}+h-\frac12)M_{\I}},& \text{J},
 \end{cases} \label{XiD}\\
 P_{\mathcal{D},n}(\eta)
 \eqdef\text{W}[\mu_{d_1},\ldots,\mu_{d_M},P_n](\eta)\times
 \begin{cases}
 \eta^{(M_{\I}+g+\frac12)M_{\II}}e^{-M_{\I}\eta},& \text{L},\\
 \bigl(\frac{1-\eta}{2}\bigr)^{(M_{\I}+g+\frac12)M_{\II}}
 \bigl(\frac{1+\eta}{2}\bigr)^{(M_{\II}+h+\frac12)M_{\I}},& \text{J},
 \end{cases} \label{PDn}\\
 \mu_{\text{v}}(\eta)\eqdef
 \begin{cases}
 e^{\eta} L^{(g-\frac12)}_{\text{v}}(-\eta),
 &\text{L, $\text{v}$ Type $\I$},\\
 \eta^{\frac12-g} L^{(\frac12-g)}_{\text{v}}(\eta),
 &\text{L, $\text{v}$ Type $\II$},\\
 \bigl(\frac{1+\eta}{2}\bigr)^{\frac12-h}
 P^{(g-\frac12,\frac12-h)}_{\text{v}}(\eta),
 &\text{J, $\text{v}$ Type $\I$},\\
 \bigl(\frac{1-\eta}{2}\bigr)^{\frac12-g}
 P^{(\frac12-g,h-\frac12)}_{\text{v}}(\eta),
 & \text{J, $\text{v}$ Type~$\II$},
 \end{cases} \label{muv}\\
 P_n(\eta) \eqdef
 \begin{cases}
 L^{(g-\frac12)}_n(\eta),& \text{L},\\
 P^{(g-\frac12,h-\frac12)}_n(\eta),& \text{J}.
 \end{cases} \label{Pn}
\end{gather}
The second is based on the Wronskians of the virtual state wavefunctions $\tilde{\phi}_{\text{v}}(x)$ and the eigenfunction $\phi_n(x)$~\cite{rrmiop3}:
\begin{gather}
\Xi_{\mathcal{D}}(\eta)=\cF^{-\frac12M(M-1)}
 \text{W}\big[\tilde{\phi}_{d_1},\ldots,\tilde{\phi}_{d_M}\big](x)\nonumber\\
 \hphantom{\Xi_{\mathcal{D}}(\eta)=}{}
 \times
 \begin{cases}
 \eta^{-M'(M'+g-\frac12)}e^{-M'\eta},& \text{L},\\
 \bigl(\frac{1-\eta}{2}\bigr)^{-M'(M'+g-\frac12)} \bigl(\frac{1+\eta}{2}\bigr)^{-M'(M'-h+\frac12)}, &\text{J},
 \end{cases} \label{W=XiD}\\
P_{\mathcal{D},n}(\eta)=\cF^{-\frac12M(M+1)}
 \text{W}\big[\tilde{\phi}_{d_1},\ldots,\tilde{\phi}_{d_M},\phi_n\big](x)\nonumber\\
\hphantom{P_{\mathcal{D},n}(\eta)=}{} \times
 \begin{cases}
 \eta^{-(M'+\frac12)(M'+g)}e^{-(M'-\frac12)\eta},& \text{L},\\
 \bigl(\frac{1-\eta}{2}\bigr)^{-(M'+\frac12)(M'+g)} \bigl(\frac{1+\eta}{2}\bigr)^{-(M'-\frac12)(M'-h)}, &\text{J},
 \end{cases} \label{W=PDn}
\end{gather}
in which $M'\eqdef\frac12(M_{\I}-M_{\II})$ and $\eta=\eta(x)$.}

\subsection[Simplif\/ied forms of the multi-indexed Laguerre and Jacobi polynomials, A]{Simplif\/ied forms of the multi-indexed Laguerre\\ and Jacobi polynomials, A} \label{sec:simpA}

In this and the subsequent subsections we will derive simplif\/ied expressions of the multi-indexed Laguerre and Jacobi polynomials $P_{\mathcal{D},n}(\eta)$ and the corresponding denominator polynomials $\Xi_\mathcal{D}(\eta)$ starting from the Wronskian expressions \eqref{XiD}--\eqref{W=PDn} in the previous subsection.

In this subsection we simplify the Wronskians of the `polynomials' \eqref{XiD}--\eqref{Pn}. Let us f\/irst transform the higher derivatives of the `seed polynomials' $\mu_{\text{v}}(\eta)$ in \eqref{muv}. For the L system, we obtain
\begin{alignat*}{3}
 & \text{L1}\colon\quad && \partial_{\eta} \bigl(e^{\eta}L^{(g-\frac12)}_{\text{v}}(-\eta)\bigr) =e^{\eta}L^{(g+1-\frac12)}_{\text{v}}(-\eta),& \\
 & \text{L2}\colon \quad & &\partial_{\eta}
 \bigl(\eta^{\frac12-g}L^{(\frac12-g)}_{\text{v}}(\eta)\bigr) =\big(\text{v}-g+\tfrac12\big)\eta^{\frac12-(g+1)} L^{(\frac12-(g+1))}_{\text{v}}(\eta),&
\end{alignat*}
by using \eqref{Lforward}--\eqref{Lid1} and \eqref{Lforward}--\eqref{Lid2}, respectively. By repeating these we arrive at for $j\in\mathbb{Z}_{\geq 1}$,
\begin{gather*}
 \partial_{\eta}^{j-1}\bigl(e^{\eta}L^{(g-\frac12)}_{\text{v}}(-\eta)\bigr) =e^{\eta} L^{(g+j-\frac32)}_{\text{v}}(-\eta),\\
\partial_{\eta}^{j-1}\bigl(\eta^{\frac12-g} L^{(\frac12-g)}_{\text{v}}(\eta)\bigr) =(-1)^{j-1}\big(g-\tfrac12-\text{v}\big)_{j-1}\eta^{\frac32-g-j}
 L^{(\frac32-g-j)}_{\text{v}}(\eta) \\
\phantom{\partial_{\eta}^{j-1}\bigl(\eta^{\frac12-g} L^{(\frac12-g)}_{\text{v}}(\eta)\bigr)}
 =\eta^{\frac12-g-K} (-1)^{j-1}\big(g-\tfrac12-\text{v}\big)_{j-1}\eta^{K+1-j} L^{(\frac32-g-j)}_{\text{v}}(\eta),
\end{gather*}
in which $K$ is a positive integer. For the J system, we obtain
\begin{alignat*}{3}
& \text{J1}\colon \quad & &\partial_{\eta}\bigl(\big(\tfrac{1+\eta}{2}\big)^{\frac12-h} P^{(g-\frac12,\frac12-h)}_{\text{v}}(\eta)\bigr)
 =\tfrac12\big(\text{v}-h+\tfrac12\big)\big(\tfrac{1+\eta}{2}\big)^{\frac12-(h+1)} P^{(g+1-\frac12,\frac12-(h+1))}_{\text{v}}(\eta),& \\
& \text{J2}\colon \quad &&\partial_{\eta}\bigl(\big(\tfrac{1-\eta}{2}\big)^{\frac12-g} P^{(\frac12-g,h-\frac12)}_{\text{v}}(\eta)\bigr)
 =-\tfrac12\big(\text{v}-g+\tfrac12\big)\big(\tfrac{1-\eta}{2}\big)^{\frac12-(g+1)} P^{(\frac12-(g+1),h+1-\frac12)}_{\text{v}}(\eta),&
\end{alignat*}
by using \eqref{Jforward},~\eqref{Jid1} and \eqref{Jforward}, \eqref{Jid1m}, respectively. Repeated applications of these formulas lead to
\begin{gather*}
 \partial_{\eta}^{j-1}\bigl(\big(\tfrac{1+\eta}{2}\big)^{\frac12-h} P^{(g-\frac12,\frac12-h)}_{\text{v}}(\eta)\bigr)
 =\frac{(-1)^{j-1}}{2^{j-1}}\big(h-\tfrac12-\text{v}\big)_{j-1} \big(\tfrac{1+\eta}{2}\big)^{\frac32-h-j} P^{(g+j-\frac32,\frac32-h-j)}_{\text{v}}(\eta)\\
\qquad{} =\big(\tfrac{1+\eta}{2}\big)^{\frac12-h-K} \frac{(-1)^{j-1}}{2^{j-1}}\big(h-\tfrac12-\text{v}\big)_{j-1}
 \big(\tfrac{1+\eta}{2}\big)^{K+1-j} P^{(g+j-\frac32,\frac32-h-j)}_{\text{v}}(\eta),\\
\partial_{\eta}^{j-1}\bigl(\big(\tfrac{1-\eta}{2}\big)^{\frac12-g} P^{(\frac12-g,h-\frac12)}_{\text{v}}(\eta)\bigr)
 =\frac{1}{2^{j-1}}\big(g-\tfrac12-\text{v}\big)_{j-1} \big(\tfrac{1-\eta}{2}\big)^{\frac32-g-j} P^{(\frac32-g-j,h+j-\frac32)}_{\text{v}}(\eta) \\
\qquad{} =\big(\tfrac{1-\eta}{2}\big)^{\frac12-g-K} \frac{1}{2^{j-1}}\big(g-\tfrac12-\text{v}\big)_{j-1} \big(\tfrac{1-\eta}{2}\big)^{K+1-j}
 P^{(\frac32-g-j,h+j-\frac32)}_{\text{v}}(\eta).
\end{gather*}
The higher derivatives of the eigenpolynomials $P_n(\eta)$ in \eqref{Pn} are replaced simply through \eqref{Lforward} and \eqref{Jforward}, respectively,
\begin{gather*}
 \partial_{\eta}^{j-1}L^{(g-\frac12)}_n(\eta)=(-1)^{j-1}L^{(g+j-\frac32)}_{n+1-j}(\eta),\\
 \partial_{\eta}^{j-1}P^{(g-\frac12,h-\frac12)}_n(\eta)=\frac{1}{2^{j-1}}(n+g+h)_{j-1} P^{(g+j-\frac32,h+j-\frac32)}_{n+1-j}(\eta),
\end{gather*}
in which we adopt the convention $L^{(\alpha)}_n(\eta)=P^{(\alpha,\beta)}_n(\eta)=0$, $n\in\mathbb{Z}_{<0}$.

Let us def\/ine $M$-dimensional column vectors $\vec{X}^{(M)}_{\text{v}}=\bigl(X^{(M)}_{\text{v},j}\bigr)_{j=1}^M$ and $\vec{Z}^{(M)}_n=\bigl(Z^{(M)}_{n,j}\bigr)_{j=1}^M$ by
\begin{gather*}
 X^{(M)}_{\text{v},j}(\eta) \eqdef
 \begin{cases}
 L^{(g+j-\frac32)}_{\text{v}}(-\eta), &\text{L, $\text{v}$ Type $\I$},\\
 (-1)^{j-1}\big(g-\tfrac12-\text{v}\big)_{j-1}\eta^{M-j} L^{(\frac32-g-j)}_{\text{v}}(\eta), &\text{L, $\text{v}$ Type $\II$},\\
 \frac{(-1)^{j-1}}{2^{j-1}}\big(h-\tfrac12-\text{v}\big)_{j-1} \big(\tfrac{1+\eta}{2}\big)^{M-j} P^{(g+j-\frac32,\frac32-h-j)}_{\text{v}}(\eta),
 &\text{J, $\text{v}$ Type $\I$},\\
 \frac{1}{2^{j-1}}\big(g-\tfrac12-\text{v}\big)_{j-1} \big(\tfrac{1-\eta}{2}\big)^{M-j} P^{(\frac32-g-j,h+j-\frac32)}_{\text{v}}(\eta),
 &\text{J, $\text{v}$ Type $\II$},
 \end{cases}\\
 Z^{(M)}_{n,j}(\eta)\eqdef
 \begin{cases}
 (-1)^{j-1}L^{(g+j-\frac32)}_{n+1-j}(\eta), &\text{L},\\
 \frac{1}{2^{j-1}}(n+g+h)_{j-1} P^{(g+j-\frac32,h+j-\frac32)}_{n+1-j}(\eta), &\text{J}.
 \end{cases}
\end{gather*}
The Wronskians in \eqref{XiD}--\eqref{PDn} are replaced by ordinary determinants consisting of these column vectors:
 \begin{gather*}
 \text{W}[\mu_{d_1},\ldots,\mu_{d_M}](\eta)= \bigl|\vec{X}^{(M)}_{d_1}(\eta) \cdots \vec{X}^{(M)}_{d_M}(\eta)\bigr|\\
\hphantom{\text{W}[\mu_{d_1},\ldots,\mu_{d_M}](\eta)= }{} \times
 \begin{cases}
 \bigl(e^{\eta}\bigr)^{M_{\I}}\bigl(\eta^{\frac32-g-M}\bigr)^{M_{\II}}, &\text{L},\\
 \bigl(\big(\frac{1+\eta}{2}\big)^{\frac32-h-M}\bigr)^{M_{\I}} \bigl(\big(\frac{1-\eta}{2}\big)^{\frac32-g-M}\bigr)^{M_{\II}}, &\text{J},
 \end{cases}\\
 \text{W}[\mu_{d_1},\ldots,\mu_{d_M},P_n](\eta)= \bigl|\vec{X}^{(M+1)}_{d_1}(\eta) \cdots \vec{X}^{(M+1)}_{d_M}(\eta)
 \vec{Z}^{(M+1)}_n(\eta)\bigr| \\
\hphantom{\text{W}[\mu_{d_1},\ldots,\mu_{d_M},P_n](\eta)=}{} \times
 \begin{cases}
 \bigl(e^{\eta}\bigr)^{M_{\I}}\bigl(\eta^{\frac12-g-M}\bigr)^{M_{\II}}, &\text{L},\\
 \bigl(\big(\frac{1+\eta}{2}\big)^{\frac12-h-M}\bigr)^{M_{\I}} \bigl(\big(\frac{1-\eta}{2}\big)^{\frac12-g-M}\bigr)^{M_{\II}}, &\text{J}.
 \end{cases}
\end{gather*}
We arrive at the main result, the simple expressions of $\Xi_{\mathcal{D}}(\eta)$ and $P_{\mathcal{D},n}(\eta)$
\begin{gather}
 \Xi_{\mathcal{D}}(\eta)=A \bigl|\vec{X}^{(M)}_{d_1}(\eta) \cdots \vec{X}^{(M)}_{d_M}(\eta)\bigr|, \label{XiDdet}\\
 P_{\mathcal{D},n}(\eta)=A \bigl|\vec{X}^{(M+1)}_{d_1}(\eta) \cdots \vec{X}^{(M+1)}_{d_M}(\eta) \vec{Z}^{(M+1)}_n(\eta)\bigr|, \label{PDndet}\\
 A= \begin{cases}
 \eta^{-M_{\II}(M_{\II}-1)}, &\text{L},\\
 \big(\frac{1+\eta}{2}\big)^{-M_{\I}(M_{\I}-1)} \big(\frac{1-\eta}{2}\big)^{-M_{\II}(M_{\II}-1)}, &\text{J}.
 \end{cases} \nonumber
\end{gather}
It should be stressed that the components of the matrices in \eqref{XiDdet} and \eqref{PDndet} are all polynomials in $\eta$. This is a good contrast
with the starting Wronskians $\text{W}[\mu_{d_1},\ldots,\mu_{d_M}](\eta)$ and $\text{W}[\mu_{d_1},\ldots,\mu_{d_M},P_n](\eta)$ in \eqref{XiD}, \eqref{PDn}, in which $\mu_{d_j}$'s have non-polynomial factors \eqref{muv}.

\subsection[Simplif\/ied forms of the multi-indexed Laguerre and Jacobi polynomials, B]{Simplif\/ied forms of the multi-indexed Laguerre\\ and Jacobi polynomials, B} \label{sec:simpB}

Here we will simplify the Wronskians of the virtual state wavefunctions and the eigenpolynomials \eqref{W=XiD}, \eqref{W=PDn}. By replacing the even order derivatives of the virtual state wavefunctions and eigenfunctions in the Wronskians \eqref{W=XiD}, \eqref{W=PDn}, by the rule $\psi^{(2m)}(x)\to(-\mathcal{E})^m\psi(x)$, we obtain
\begin{gather}
 \text{W}\big[\tilde{\phi}_{d_1},\ldots,\tilde{\phi}_{d_M},\phi_n\big](x) =\det(a_{j,k})_{1\leq j,k\leq M+1}, \label{W1}\\
 a_{2l-1,k}=(-\tilde{\mathcal{E}}_{d_k})^{l-1}\tilde{\phi}_{d_k}(x),\qquad 1\leq k\leq M,\qquad 1\leq l\leq\big[\tfrac{M+2}{2}\big],\nonumber\\
a_{2l-1,M+1}=(-\mathcal{E}_n)^{l-1}\phi_n(x), \qquad 1\leq l\leq\big[\tfrac{M+2}{2}\big],\nonumber\\
 a_{2l,k} =(-\tilde{\mathcal{E}}_{d_k})^{l-1}\tilde{\phi}'_{d_k}(x), \qquad 1\leq k\leq M,\qquad 1\leq l\leq\big[\tfrac{M+1}{2}\big],\nonumber\\
a_{2l,M+1}=(-\mathcal{E}_n)^{l-1}\phi'_n(x), \qquad 1\leq l\leq\big[\tfrac{M+1}{2}\big], \nonumber
\end{gather}
in which $[a]$ denotes the greatest integer not exceeding $a$. The f\/irst derivatives in the $2l$-th row can be simplif\/ied by adding $-\frac{\phi_0'(x)}{\phi_0(x)}\times(2l-1)$-th row,
\begin{gather}
 \phi'_n(x)\to\left(\frac{d}{dx}-\frac{\phi'_0(x)}{\phi_0(x)}\right)\phi_n(x)=\frac{\cF}{\eta'(x)}\phi_0(x)\zeta_n\bigl(\eta(x)\bigr)
 =\phi_0(x)\zeta_n\bigl(\eta(x)\bigr) A,\nonumber\\
 \tilde{\phi}'_{\text{v}}(x)\to \left(\frac{d}{dx}-\frac{\phi'_0(x)}{\phi_0(x)}\right) \tilde{\phi}_{\text{v}}(x)
 =\frac{\cF}{\eta'(x)}\,\tilde{\phi}_0(x) \tilde{\zeta}_{\text{v}}\bigl(\eta(x)\bigr) =\tilde{\phi}_0(x)\tilde{\zeta}_{\text{v}}\bigl(\eta(x)\bigr) A,\nonumber\\
 A=\frac{\cF}{\eta'(x)}=
 \begin{cases}
 x^{-1}, & \text{L},\\
 (\sin x\cos x)^{-1}, & \text{J},
 \end{cases} \label{phi'}
\end{gather}
in which $\zeta_n(\eta)$ and $\tilde{\zeta}_{\text v}(\eta)$ are polynomials in $\eta$ def\/ined by
\begin{gather*}
 \zeta_n(\eta)\eqdef
 \begin{cases}
 -2\eta L^{(g+\frac12)}_{n-1}(\eta), &\text{L},\\
 -\tfrac12(n+g+h)\big(1-\eta^2\big)P^{(g+\frac12,h+\frac12)}_{n-1}(\eta), &\text{J},
 \end{cases}\\
 \tilde{\zeta}_{\text{v}}(\eta) \eqdef
 \begin{cases}
 2\eta L^{(g+\frac12)}_{\text{v}}(-\eta), &\text{L, $\I$},\\
 -2\big(g-\tfrac12-\text{v}\big)L^{(-g-\frac12)}_{\text{v}}(\eta), &\text{L, $\II$},\\
 \big(h-\tfrac12-\text{v}\big)(1-\eta) P^{(g+\frac12,-h-\frac12)}_{\text{v}}(\eta), &\text{J, $\I$},\\
 -\big(g-\tfrac12-\text{v}\big)(1+\eta) P^{(-g-\frac12,h+\frac12)}_{\text{v}}(\eta), & \text{J, $\II$}.
 \end{cases}
\end{gather*}
Use is made of \eqref{Lforward}--\eqref{Lid2} for L and \eqref{Jforward}--\eqref{Jid1m} for J.

By extracting the functions $\phi_0(x)$, $\tilde{\phi}_0(x)$ from each column of the matrix $a_{j,k}$ \eqref{W1} and the factor $A$ of \eqref{phi'} from
the even rows, we arrive at another set of simplif\/ied determinant expressions for the multi-indexed polynomials:
\begin{gather}
P_{\mathcal{D},n}(\eta)=\cF^{-\frac12M(M+1)} \det(a_{j,k})_{1\leq j,k\leq M+1} A, \label{PDn1}\\
 a_{2l-1,k}=(-\tilde{\mathcal{E}}_{d_k})^{l-1}\xi_{d_k}(\eta), \qquad 1\leq k\leq M, \qquad 1\leq l\leq\big[\tfrac{M+2}{2}\big],\nonumber\\
a_{2l-1,M+1} =(-\mathcal{E}_n)^{l-1}P_n(\eta), \qquad 1\leq l\leq\big[\tfrac{M+2}{2}\big],\nonumber\\
 a_{2l,k} =(-\tilde{\mathcal{E}}_{d_k})^{l-1}\tilde{\zeta}_{d_k}(\eta), \qquad 1\leq k\leq M,\qquad 1\leq l\leq\big[\tfrac{M+1}{2}\big],\nonumber\\
a_{2l,M+1} =(-\mathcal{E}_n)^{l-1}\zeta_n(\eta), \qquad 1\leq l\leq\big[\tfrac{M+1}{2}\big],\nonumber\\
 A= \begin{cases}
 \eta^{-([M']+1)([M']+M-2[\frac{M}{2}])}, &\text{L},\\
 \bigl(\tfrac{1-\eta}{2}\bigr)^{-([M']+1)([M']+M-2[\frac{M}{2}])} \bigl(\tfrac{1+\eta}{2}\bigr)^{-([-M']+1)([-M']+M-2[\frac{M}{2}])}, &\text{J}.
 \end{cases} \nonumber
\end{gather}
In a similar way, we obtain
\begin{gather}
 \Xi_{\mathcal{D}}(\eta)=\cF^{-\frac12M(M-1)} \det(a_{j,k})_{1\leq j,k\leq M} A, \label{XiD1}\\
 a_{2l-1,k}=(-\tilde{\mathcal{E}}_{d_k})^{l-1}\xi_{d_k}(\eta), \qquad 1\leq l\leq\big[\tfrac{M+1}{2}\big], \nonumber\\
 a_{2l,k} =(-\tilde{\mathcal{E}}_{d_k})^{l-1}\tilde{\zeta}_{d_k}(\eta), \qquad 1\leq l\leq\big[\tfrac{M}{2}\big], \nonumber\\
 A=
 \begin{cases}
 \eta^{-[M']([M']+M-2[\frac{M}{2}])}, & \text{L},\\
 \bigl(\tfrac{1-\eta}{2}\tfrac{1+\eta}{2}\bigr)^{-[M']([M']+M-2[\frac{M}{2}])}, &\text{J}.
 \end{cases} \nonumber
\end{gather}
Again all the components of the matrices $a_{j,k}$ in \eqref{PDn1} and \eqref{XiD1} are polynomials in $\eta$.

\subsection{Parity transformation of the multi-indexed Jacobi polynomials}\label{sec:pari}

The Jacobi polynomial has the parity transformation property \cite{szego}
\begin{gather}
 P^{(\alpha,\beta)}_n(-x)=(-1)^nP^{(\beta,\alpha)}_n(x). \label{P:J}
\end{gather}
We will show that this property is inherited by the multi-indexed Jacobi polynomials. It is based on the property of the Wronskian
\begin{gather*}
 \text{W}[f_1,\ldots,f_n](-\eta) =(-1)^{\frac12n(n-1)}\text{W}[g_1,\ldots,g_n](\eta),\qquad g_k(\eta)\eqdef f_k(-\eta).
\end{gather*}

In this subsection, we indicate the types of the virtual states explicitly by $(\text{v},\text{t})$, in which $\text{t}$ stands for Type~$\I$ or~$\II$. Based on \eqref{P:J}, we obtain
\begin{gather*}
 \mu_{(\text{v},\text{t})}\bigl({-}\eta;(g,h)\bigr) =(-1)^{\text{v}}\mu_{(\text{v},\bar{\text{t}})}\bigl(\eta;(h,g)\bigr),\qquad
 \bar{\text{t}}\eqdef
 \begin{cases}
 \II, & \text{t}=\I,\\
 \I, &\text{t}=\II,
 \end{cases}\\
 P_n\bigl({-}\eta;(g,h)\bigr) =(-1)^nP_n\bigl(\eta;(h,g)\bigr).
\end{gather*}
For the multi-index set $\mathcal{D}=\{(d_1,\text{t}_1),\ldots,(d_M,\text{t}_M)\}$ of the virtual state wavefunctions, let us def\/ine the `mirror ref\/lected' multi-index set $\mathcal{D}'\eqdef\{(d_1,\bar{\text{t}}_1),\ldots,(d_M,\bar{\text{t}}_M)\}$. Corresponding to $M_{\I}\eqdef\#\{d_j\, |\, (d_j,\I)\in\mathcal{D}\}$, $M_{\II}\eqdef\#\{d_j\,|\,(d_j,\II)\in\mathcal{D}\}$, we have $M'_{\I}\eqdef\#\{d_j\,|\,(d_j,\I)\in\mathcal{D}'\}=M_{\II}$,
$M'_{\II}\eqdef\#\{d_j\,|\, (d_j,\II)\in\mathcal{D}'\}=M_{\I}$. By parity transformation $\eta\to-\eta$, the multi-indexed Jacobi polynomial $P_{\mathcal{D},n}\bigl(\eta;(g,h)\bigr)$ is mapped to $P_{\mathcal{D}',n}\bigl(\eta;(h,g)\bigr)$ with the `mirror ref\/lected' multi-index set~$\mathcal{D}'$:
\begin{gather*}
 P_{\mathcal{D},n}\bigl({-}\eta;(g,h)\bigr) \\
 \qquad{}=\big(\tfrac{1-\eta}{2}\big)^{(M_{\II}+h+\frac12)M_{\I}}
 \big(\tfrac{1+\eta}{2}\big)^{(M_{\I}+g+\frac12)M_{\II}} \text{W}\big[\mu_{(d_1,\text{t}_1)},\ldots,\mu_{(d_M,\text{t}_M)},P_n\big]
 \bigl({-}\eta;(g,h)\bigr)\nonumber\\
 \qquad{} =\big(\tfrac{1-\eta}{2}\big)^{(M_{\II}+h+\frac12)M_{\I}}
 \big(\tfrac{1+\eta}{2}\big)^{(M_{\I}+g+\frac12)M_{\II}} (-1)^{\frac12M(M+1)}\nonumber\\
 \qquad\quad{} \times
 \text{W}\big[(-1)^{d_1}\mu_{(d_1,\bar{\text{t}}_1)},\ldots, (-1)^{d_M}\mu_{(d_M,\bar{\text{t}}_M)},(-1)^nP_n\big]\bigl(\eta;(h,g)\bigr)\nonumber\\
 \qquad{}=\big(\tfrac{1-\eta}{2}\big)^{(M_{\II}+h+\frac12)M_{\I}} \big(\tfrac{1+\eta}{2}\big)^{(M_{\I}+g+\frac12)M_{\II}} (-1)^{\frac12M(M+1)}(-1)^{d_1+\cdots+d_M+n} \\
 \qquad\quad{}\times \text{W}\big[\mu_{(d_1,\bar{\text{t}}_1)},\ldots,\mu_{(d_M,\bar{\text{t}}_M)}, P_n\big]\bigl(\eta;(h,g)\bigr) \\
 \qquad{} =(-1)^{n+\frac12M(M+1)+\sum\limits_{k=1}^Md_k} P_{\mathcal{D}',n}\bigl(\eta;(h,g)\bigr).
\end{gather*}
Similarly, the denominator polynomial $\Xi_{\mathcal{D}}\bigl(\eta;(g,h)\bigr)$ is mapped to $\Xi_{\mathcal{D}'}\bigl(\eta;(h,g)\bigr)$ with a sign factor:
\begin{gather*}
 \Xi_{\mathcal{D}}\bigl({-}\eta;(g,h)\bigr) =(-1)^{\frac12M(M-1)+\sum\limits_{k=1}^Md_k} \Xi_{\mathcal{D}'}\bigl(\eta;(h,g)\bigr).
\end{gather*}
For the special case of `mirror symmetric' multi-index set $\mathcal{D}'=\mathcal{D}$ (as a set), i.e., $\{d_j\,|\,(d_j,\I)\in\mathcal{D}\}= \{d_j\,|\,(d_j,\II)\in\mathcal{D}\}$ (as a set), we have $P_{\mathcal{D}',n}(\eta)=\pm P_{\mathcal{D},n}(\eta)$. In fact, this formula turns out to be
\begin{gather*}
 P_{\mathcal{D}',n}(\eta)=(-1)^{(\frac{M}{2})^2}P_{\mathcal{D},n}(\eta).
\end{gather*}
For this special case the parity transformation gives
\begin{gather*}
 P_{\mathcal{D},n}\bigl({-}\eta;(g,h)\bigr) =(-1)^nP_{\mathcal{D},n}\bigl(\eta;(h,g)\bigr),\qquad \Xi_{\mathcal{D}}\bigl({-}\eta;(g,h)\bigr)
 =\Xi_{\mathcal{D}}\bigl(\eta;(h,g)\bigr).
\end{gather*}

\section{Summary and comments}\label{sec:summary}

The multi-indexed Laguerre and Jacobi polynomials are def\/ined by the Wronskian expressions originating from multiple Darboux transformations. Two simplif\/ied determinant expressions of them, \eqref{XiDdet},~\eqref{PDndet} and \eqref{PDn1},~\eqref{XiD1}, which do not contain derivatives, are derived based on the pro\-perties of the Wronskian and identities of the Laguerre and Jacobi polynomials. For~\eqref{PDn1},~\eqref{XiD1}, the Schr\"odinger equation is used. For~\eqref{XiDdet},~\eqref{PDndet}, various identities of the Laguerre and Jacobi polynomials are used, which are essentially forward shift relations. Although the calculation in Section~\ref{sec:simpA} is performed for polynomials, it can be done for wavefunctions just like~\cite{os29}, in which simplif\/ied determinant expressions are presented for the multi-indexed polynomials obtained by multiple Darboux transformations with pseudo virtual states wavefunctions as seed solutions. The parity transformation property of the multi-indexed Jacobi polynomials is also derived.

Multi-indexed orthogonal polynomials have been constructed for the classical orthogonal polynomials in the Askey scheme \cite{askey,kls}, i.e., the Wilson, Askey--Wilson, Meixner, little $q$-Jacobi and ($q$-)Racah polynomials, etc.~\cite{os26,os27,os35}. These polynomials belong to `discrete' quantum mechanics~\cite{os24}, in which the Schr\"odinger equations are second-order dif\/ference equations. The Casoratian expressions of these multi-indexed polynomials can also be simplif\/ied by using various identities as demonstrated here. These simplif\/ications will be published elsewhere \cite{odake16}.

\subsection*{Acknowledgements}

S.O.~is supported in part by Grant-in-Aid for Scientif\/ic Research from the Ministry of Education, Culture, Sports, Science and Technology (MEXT), No.~25400395.

\pdfbookmark[1]{References}{ref}
\LastPageEnding

\end{document}